\documentclass[10pt]{amsart}
\usepackage{color}
\usepackage{hyperref}
\usepackage{amsmath,amsfonts,latexsym,amssymb}
\usepackage[latin1]{inputenc}
\usepackage{amsmath}
\usepackage{graphicx}
\newtheorem{theorem}{Theorem}[subsection]
\newtheorem{thm}{Theorem}[subsection]
\newtheorem{lemma}[theorem]{Lemma}

\newtheorem{coro}[theorem]{Corollary}

\def\bb{\begin{eqnarray*}}
\def\ee{\end{eqnarray*}}

\def\ZZ{\mathbb Z}
\def\HH{\Bbb H}

\def\MM{\Sigma^*}

\def\tr{\mathrm{tr}\,}

\def\RR{\mathbb{R}}

\def\ff{f^\sim}

\def\tr{\mathrm{tr}\,}

\def\sl2{SL_{2}(\RR)}

\def\ra{\rightarrow}

\def\dhyp{\partial \HH}

\def\fund{\pi_{1}}

\def\inv{I}

\begin{document}

%
%

\title[Normalized Entropy versus Volume]
{Large cone angles on a punctured sphere.}

\author[G. McShane]{%
	Greg McShane 
} 
\address{%
	UFR de Math\'ematiques \\
	Institut Fourier 100 rue des maths \\
	BP 74, 38402 St Martin d'H\`eres cedex, France
}
\email{%
	Greg.McShane@ujf-grenoble.fr
}


\thanks{%
The author is partially supported by ANR MODGROUP
} 

\begin{abstract} 
Do and Norbury found  a so-called differential relation 
which relates the volume of  the moduli space of  singular surface with a cone point
to that of a smooth surface obtained by forgetting the cone point.
Their procedure is valid for cone angles less  than $\pi$ by work of Tan, Wong and Zhang.
We study the moduli space of a surface with a single cone point of angle 
ranging from $0$ to $2\pi$ 
using a coordinate system closely related to Penner's $\lambda$-lengths.
We compute the action of the mapping class group in these coordinates,
give an explicit expression for Wolpert's symplectic form
and use this to justify Do and Norbury's approach using 
just hyperbolic geometry.

\end{abstract}

\maketitle

\section{Introduction}

The aim of this paper is to extend the results of Nakanishi and Naatanen \cite{NaNa}
to hyperbolic cone surfaces of signature
$(0;  0 , 0 , 0 ,\theta)$ where $\theta$ is a cone angle in the interval $[\pi,2\pi[$.
Our motivation for this is to understand geometrically the "differential recurrences"
introduced by Do and Norbury \cite{DoNo}.
In particularly we want to see how the hyperbolic structure degenerates 
as the cone angle $\theta$ approaches $2\pi$.

We use a different coordinate system to Nakanishi and Naatanen,
who use traces of matrices,
and Maloni, Palesi and Tan \cite{MaPaTa}  also study the 4-holed sphere 
using similar ideas.
Each of our coordinates is a cross ratio and so our results should
generalise to representations  into higher rank groups 
as in \cite{McLab} and \cite{Tan}.
The price to pay for this (see Section \ref{autos})
is that determining the action of the (extended) mapping class group
is slightly more delicate in our coordinates.
In the trace coordinates the generators are obtained via
so-called root flipping or Vieta jumping.
Finally we note   that \cite{SchTra} have a different approach,
based on generalising quadratic differentials,
 to the problems posed by large cone angles.

The main difficulty is that, although it is easy to construct a hyperbolic cone surfaces with this signature,
it is not clear that given any other representation in the same  component of the relative character variety
that it is also the holonomy of a (unique) hyperbolic structure
 (Theorem \ref{metric}).

\subsection{Surfaces with  geodesic boundary or cone points}

As in \cite{tanzhang}, \cite{DoNo} we consider a cone point on a hyperbolic surface $\MM$
to be a ``generalized  boundary component" having purely imaginary length.

\begin{enumerate}
	\item If  the surface  $\MM$ has  a single puncture and $[\delta] \subset  \fund(\MM)$
 represents a the boundary of a small embedded disc around the puncture
 then  its Teichmuller space embeds as a connected subset of the $\sl2$-character variety of $\fund$ 
containing representations such that  $|\tr \rho(\delta)| = 2$.
By convention $\delta$ represents a boundary component of zero length. 
	\item If   $\MM$ has  a single totally geodesic boundary component  
represented by a  conjugacy class of loops $[\delta] \subset  \fund(\MM)$,
the space of marked hyperbolic structures on $\MM$
 such that the boundary geodesic has length $\ell_\delta >0 $
embeds in the $\sl2$-character variety of its fundamental group.
The image is a connected component of the \textit{ relative $\sl2$-character variety of $\fund$},
that is, the subvariety such that $|\tr \rho(\delta)| = 2 \cosh(\ell_\delta)$.
	\item If the puncture is ``replaced" by  a single cone point of angle $\theta \leq \pi$ \cite{tanzhang},
one  identifies the space  of hyperbolic structures on the  surface $\MM$ 
with a component 
of the $\sl2$-character variety 
such that $|\tr \rho(\delta)| = 2 \cos(\theta/2)$.
It proves  useful to adopt the convention that $\delta$
represents a boundary component of purely imaginary length  $\ell_\delta = i \theta$.
\end{enumerate}
In each of these  three cases we say that $\delta$ is a \textit{generalized boundary component.}

We consider representations $\rho$ obtained by deforming 
the holonomy of a metric on a surface  with a single cone point of angle $\theta \leq 2\pi$
through representations such that $|\tr \rho(\delta)| = |2 \cos(\theta/2)|$.
The set of such representation forms a subset of the $\sl2$-character variety
which we call a \textit{geometric component of the  relative character variety}.
Note that, if $\theta > \pi$ then  none of the  representation 
$\rho$ of $\fund$  is  discrete and faithful.
Consequently, we  may not assume that
given an essential simple loop $\gamma$  
and a representation $\rho$ the isometry $\rho(\gamma)\in \sl2$ will be hyperbolic.

With this in mind, our aim is to show that  the geometric component of the  relative character variety
enjoys the following properties which should be familiar from the theory of Teichmuller space;

\begin{enumerate}
	\item every simple loop is hyperbolic.
	\item the action of the outer automorphism group of $\fund$ is proper. 
	\item the set of lengths of simple loops has Fibonnaci growth.
	\item every representation is the holonomy of a hyperbolic structure with a single cone point on a surface.
\end{enumerate}

\subsection{Generalised four holed spheres}

The fundamental group  of the four holed  sphere is  freely generated 
by three  peripheral loops
$\alpha,\beta,\gamma$.
After possibly replacing some  of these elements
by their inverses,
 we may assume that
$\delta = \gamma\beta\alpha$ is also a peripheral element.
Throughout each of  $\alpha,\beta,\gamma$ will be either hyperbolic or parabolic element of $\sl2$
with positive trace. 
On the other hand we impose no restriction on the type of  $\delta$ 
but its trace will always be strictly greater than $-2$.

We make the following definition:  a  component 
(of the relative character variety $|\tr \rho(\delta)| = $ constant) 
is \textit{geometric} iff it contains the holonomy of a metric of  a cone surface.
For example, if one takes  $a= b= c = 3$ then $\kappa(3,3,3) = 2$,
the matrices $A,B,C$ are  in the congruence subgroup $\Gamma_3 < SL(2,\ZZ)$
and $\HH/ \langle A,B,C\rangle$ is the regular 4 punctured sphere.

\begin{theorem}\label{thmA}
For the geometric component of the character variety of the generalised 4 holed sphere

\begin{enumerate}
	\item every simple loop is hyperbolic.
	\item the action of the outer automorphism group is proper.
	\item the set of lengths of simple loops has Fibonnaci growth.
	\item every representation is the holonomy of a hyperbolic structure with a single cone point on a surface.
\end{enumerate}
\end{theorem}

\subsection{An inequality for generalised four holed spheres}

If the four holed  sphere is equipped with a complete hyperbolic metric, 
and the boundary $\delta$ is totally geodesic of lengths $\ell_\delta >0$ 
then the loops $\alpha\beta$ and $\beta\gamma$ can be represented 
by closed simple geodesics which meet in exactly two points and which satisfy
\begin{equation}\label{mmc2}
\sinh(\ell_{\alpha\beta}/4)\sinh(\ell_{\beta\gamma}/4) \geq 1.
\end{equation}
This inequality follows almost immediately from the Collar Lemma 
but we will deduce it from a more general result which we explain now.

As  discussed in the preceding paragraph,
we view  a three holed sphere with a single cone point of angle $\theta <2\pi$
as a four holed sphere with a  generalized boundary component, represented by an element
 $\delta \in \fund$ such that $\tr\rho(\delta) = 2\cos(\theta/2).$
Under the hypothesis that $\theta < 2\pi$,
the elements $\alpha\beta$ and $\beta\gamma$ are always hyperbolic (Theorem \ref{thmA}) 
so that there are real numbers
$\ell_{\alpha\beta}, \ell_{\beta\gamma}$ such that 
$$\tr \rho(\alpha\beta)  =  -2\cosh(\ell_{\alpha\beta}/2),\,\tr \rho(\beta\gamma)  =  -2\cosh(\ell_{\beta\gamma}/2)\ $$
and we prove 

\begin{theorem}\label{thmB}
With the above notation we have the following inequality
\begin{equation}\label{mmcd}	
\sinh(\ell_{\alpha\beta}/4)\sinh(\ell_{\beta\gamma}/4) 
\geq \cos(\theta/4).
\end{equation}
\end{theorem}
It is worth taking a  moment to consider the extremal cases $\theta =0,2\pi$:
\begin{itemize}
\item when  $\theta = 0$ this corresponds to a cone angle of zero,
 in other words a cusp, and we recover the inequality (\ref{mmc2}) above.
\item when  $\theta = 2\pi$ one has 
$\sinh(\ell_{\alpha\beta}/4)\sinh(\ell_{\beta\gamma}/4) \geq \cos(\pi/2) = 0. $
so that, as one expects, the inequality is redundant.
\end{itemize}

\subsection{Motivation: volume computations}

The essential motivation for this paper is understanding Do and Norbury's approach  to volume recurrences and Maloni, Palesi and Tan  \cite{MaPaTa}
study the mapping class group action on the  four holed sphere.

Nakanishi and Naatanen determined the
 symplectic volume of the moduli spaces of the 
 once holed torus $V_{1}(l_1)$
 and the 4 holed sphere $V_0(l_1,l_2,l_3,l_4)$.
 The volumes are quadratic polynomials in the lengths 
 of the boundary components :
 $$V_{1}(l_1) = {1 \over 24} (4\pi^2 + l_1^2) ,\, \, 
 V_0(l_1,l_2,l_3,l_4) = {1 \over 2}(4 \pi^2 + \sum_i l_i^2).$$
Their method consists of finding a fundamental region for the action of the mapping class group
on the relative character variety when the peripheral loops are hyperbolic or possibly parabolic.
They do this by  imitating Wolpert's \cite{wolp} analysis of the mapping class group action on the
Teichmuller space of the once punctured torus.

Do and Norbury found  a so-called differential relation 
which, in the simplest case, relates
the volume of  the moduli space of  singular surface with a cone point
to that of a smooth surface obtained by forgetting the cone point.
For a 3 punctured  sphere with a cone point of angle $\theta$
 the volume of the moduli space is
 $$V_0(0,0,0, i\theta ) = {1 \over 2}(4 \pi^2  - \theta^2),$$
 and we see that for $\theta = 2\pi$ the value of $V_0(0,0,0, 2\pi i  )$
 is $0$ which
 is  the volume of the moduli space of  the the 3-punctured  sphere
 is $0$.
 
 A more interesting case is that of the punctured torus
 with a single cone point of angle  $\theta$,
 the volume, computed by Mirzakhani is:

$$V_1(i\theta,l_2) =  {1 \over 192}(4 \pi^2 - \theta^2 + l_2^2)(12 \pi^2 -\theta^2 + l_2^2)  $$
 so for $\theta = 2\pi$ the value of $V_1(i\theta,0) $
 is a quartic in $l_2$
 is  \textit{not} the volume of the moduli space of  
 the one holed  torus.

However, if we take the derivative at $2\pi i$ we will recover 
the volume of the one holed torus as a factor:
\begin{eqnarray*}
V_1(l_1,l_2) &=& {1 \over 192}(4 \pi^2 + l_1^2 + l_2^2)(12 \pi^2 + l_1^2 + l_2^2)  \\
{d \over dl_1} V_1(l_1,l_2) &=& { l_1  \over 96} (16  \pi^2 + 2 l_1^2  + 2l_2^2)\\
\left.  {d \over dl_1}\right|_{2\pi i} V_1(l_1,l_2) &=& { 2\pi i  \over 96} (8 \pi^2 + + 2l_2^2)\\
&=& { 2\pi i  \over 4.24} (4\pi^2 + l_2^2)\\
&=& { 2\pi i  \over 4} V_1(l_1)
\end{eqnarray*}

As mentioned above Maloni, Palesi and Tan  \cite{MaPaTa}
study the mapping class group action on the 
 the (relative) $SL(2,\mathbb{C})$ character varieties of the four-holed sphere.
They  describe a domain of discontinuity, and, in the case of real characters, 
show that this domain of discontinuity may be non-empty 
on the components where the relative euler class is non-maximal.
This is to be expected, at least heuristically,
as the volume polynomial ${1 \over 2}(4 \pi^2 + \sum_i l_i^2)$ does not vanish 
on these non-maximal components.
They  should conjecturally decompose as a union of wandering domains
and a set where the mapping class action is ergodic.




\section{Moebius maps and fixed points}

Our methods are based on elementary algebra and geometry.
It is important, however, to recall some basic facts and definitions.

Recall that if $M = \begin{pmatrix}
M_{1,1}  & M_{1,2}  \\
M_{2,1} & M_{2,2}
\end{pmatrix}  \in \sl2$
then it acts on $\HH$ by mobius transformation 
$$ z \mapsto 
{
M_{1,1}z + M_{1,2} \over
M_{2,1} z + M_{2,2}}.$$
The  fixed points of Moebius transformation
are given by the following formula:

\begin{eqnarray*}
z_\pm &=& { (M_{1,1} - M_{2,2}) \pm \sqrt{ \tr^2 M - 4 } \over M_{2,1} } \\
&=& { (\tr M  - 2M_{2,2}) \pm \sqrt{ \tr^2 M - 4 } \over M_{2,1} }
\end{eqnarray*}
Note further that, if $M$ is {\it elliptic} , that is $|\tr M | <2$, then $z_\pm$ are complex conjugate and
\begin{equation}\label{fpsign}
\text{Re} (z_\pm)  
=  { \tr M  - 2M_{2,2}  \over M_{2,1} } 
\end{equation} 
on the othe hand  if $M$ is {\it hyperbolic},  $|\tr M | > 2$,  then the fixed points are points of the extended real line
$\RR \cup \{ \infty \}$.

\section{The generalized four punctured sphere}

Following the convention,
one views a 3 punctured  sphere with a single cone point 
as a \textit{generalized 4 holed sphere} 
that is a four holed sphere with three boundary components of length zero
and  a single  generalized boundary component $\delta$ of purely imaginary length.

The fundamental group $\fund$ of the 4 punctured sphere is  freely generated by 3 peripheral loops
$\alpha,\beta,\gamma$.
After possibly replacing some  of these elements
by their inverses,
 we may assume that
$\delta = \gamma\beta\alpha$ is also a peripheral element 
and we have the presentation

$$\fund = \langle \alpha,\beta,\gamma, \delta |
  \delta = \gamma\beta\alpha\rangle.$$

 \subsection{Explicit matrices and a parametrization}\label{expMatrices}

We begin by studying the set of representations of $\fund = \langle \alpha,\beta,\gamma, \delta |
  \delta = \gamma\beta\alpha\rangle$.
The group is freely generated by $\alpha,\beta,\gamma$
 and we may  define a representation
$\rho: \fund \rightarrow SL(2,\RR)$ by setting
$$\rho(\alpha) = A,\rho(\beta) = B,\rho(\gamma) = C$$
with $ A, B, C \in SL(2,\RR)$ parabolic with distinct fixed points.
Recall that an element of $SL(2,\RR)$ is \textit{ parabolic} iff it is conjugate to
$$\pm \begin{pmatrix}
1 & p \\
0 & 1
\end{pmatrix},\,p\neq 0
$$
and that such a transformation has a single (ideal) fixed point in $\dhyp$.
The action of $\sl2$ is transitive on triples of distinct points 
in $\dhyp$
and we normalize so that  the fixed points of $A,B,C$ are respectively fixes $0,1,\infty$.
So that, for some $a,b,c>0$
$$C= \begin{pmatrix}
1 & -c \\
0 &  1
\end{pmatrix},
B= \begin{pmatrix}
1+b & -b \\
b &  1-b
\end{pmatrix},
A= \begin{pmatrix}
1 & 0 \\
a&  1
\end{pmatrix}.
$$
Since $\fund$ is freely generated by $\alpha,\beta,\gamma$,
the  triple $(A,B,C)$ determines a point of the $\sl2$ representation variety of $\fund$
and we denote $[(A,B,C)]$ the corresponding point in the character variety.
Thus we have a map from $\RR^3$ to the  representation variety
$$(a,b,c) \mapsto (A,B,C) .$$ 
and from $\RR^3$ to the character variety
$$(a,b,c)  \mapsto [(A,B,C)].$$
The traces of elements of $\rho(\pi_1)$ form a natural class of
functions on the representation  variety and
one can compute explicit expressions for these functions  in terms of $a,b,c$. 
For example 
$$
BA= \begin{pmatrix}
1+b -ab  & -b \\
b+a -ab &  1-b
\end{pmatrix},$$
so that 
$$\tr AB = \tr BA = 2- ab.$$
Likewise, by explicit computation, one has
$$
\tr AC =  2- ac ,\,\,
\tr BC =  2- bc.
$$
Recall that an element of $\sl2$ is hyperbolic iff $|\tr| >2$.
For all our representations  $\tr A = \tr B = 2 >0$,
so whenever $\langle A,B\rangle$ is discrete, 
by the Triple Trace Theorem \cite{maskit}, 
\begin{equation}\label{TTT}
2 - ab = \tr AB \leq -2.
\end{equation}
We shall consider non discrete representations $\rho$
but we will show (Lemma \ref{product}) that the inequality (\ref{TTT}) 
still holds  for all the representations considered provided 
$\tr CBA > - 2$.

Finally we check that our map gives a parameterization of the character variety:

\begin{lemma}[Smooth parameterization]
The restriction of the map
\begin{eqnarray*}
(a,b,c) &\mapsto&  [(A,B,C)] \\
\end{eqnarray*}
to $\{ab, bc, ac >2 \}$ is a diffeomorphism onto its image.

\end{lemma}
\proof It suffices to check that the map
$$(a,b,c) \mapsto  
(x,y,z) = (\tr BC,\tr AC, \tr AB)
= ( 2 - bc, 2 - ac, 2 - ab)$$
is a diffeomorphism.
Both the jacobian and the inverse of this map are easy to compute and we leave this to the reader to check. 
$\Box$

\subsection{Geometric interpretation of $a,b,c$}\label{geomInterp}

If the three  matrices $A,B,C$ generate a fuchsian group $\Gamma$
then the ideal triangle $\ 0,1,\infty$ embeds in the quotient surface $\HH / \Gamma$.
The surface $\HH / \Gamma$  has  3 (and possibly 4 if $CBA$ is parabolic) cusps, 
that is  one for each of the ideal vertices of $\ 0,1,\infty$.
Each of the the edges of the triangle $\ 0,1,\infty$ 
embeds as an arc joining distinct pairs of  cusps
and each of these cusps, by Shimura's Lemma, 
lies  in one of 3 pairwise disjoint cusp regions of area 1.
In this paragraph we give an interpretation of $a,b,c$ in terms of this configuration.

The region $\{ \mathrm{Re} z \geq  c \} \subset \HH$ embeds as a cusp region of area 
$\int_{-c}^0 { {dx} \over   c} = 1$
and the portion of the  triangle $\ 0,1,\infty$ contained in this region is 
$ \int_0^1 {{dx} \over c} =  {1 \over c}.$
Such a portion, 
that is the intersection of an ideal triangle with a horoball
centered at one of its vertices and meeting exactly two of the  sides of the triangle,
is often called a \textit{prong}.
By a similar argument, one sees that the prong of the  triangle $\ 0,1,\infty$
contained in the lift of a cusp region based at $0$ (resp. $1$) 
is exactly $ {1 \over a}$ (resp. $ {1 \over b}$).

On each side of an ideal triangle there is a well-defined midpoint.
A simple calculation yields that 
the  distance from the midpoint to the prong of area $h$ is $\log(h)$.
It follows that the arc, that is the portion of the geodesic $0,\infty$,
outside of lift of the cusp region of area $1/c$ based at $\infty$
and the cusp region of area $1/a$ based at $0$,
is $\log(ac)$. In this way one sees that the quantity $ac$ is 
one of Penner's  $\lambda$-length (see \cite{penner}).

\subsection{Topology of the relative character variety}

A relative character variety is defined to be a level set of 
the function $\kappa$
where
$$\kappa(a,b,c) := \tr  CBA.$$
In general this has more than one connected component
(of codimension $0$ 
and  possibly infinitely many of codimension $>1$).
Here we are interested in the codimension $0$ component 
parameterized by  $(\RR^+)^3$ under the map
$(a,b,c) \mapsto [(A,B,C)]$ and for this parameterization
$$\kappa(a,b,c) =  2  + abc -ab - bc - ac.$$

Here we are interested a codimension $0$ component 
which contains at least one representation which is the  holonomy of  metric on a cone surface.

\noindent
\textbf{Definition} A component (of the character variety) 
is \textit{geometric} if it contains the holonomy of a metric
on a cone surface.\\

\begin{description}

\item[Example 1:] with $a= b= c = 3,\, \kappa(3,3,3) = 2$, the matrices $A,B,C$ are  in
the congruence subgroup $\Gamma_3 < SL(2,\ZZ)$
and $\HH/ \langle A,B,C\rangle$ is the regular 4 punctured sphere.

\item[Example 2:] with $a= b= c = 2,\, \kappa(2,2,2) = -2$, the matrices $A,B,C$ are  in
the congruence subgroup $\Gamma_2 < SL(2,\ZZ)$
and $\HH/ \langle A,B,C\rangle$ is the 3 punctured sphere.
Note that the image is not a free group since 
$CBA$ is a torsion of order 2.

\end{description}

\begin{lemma}[Subvariety]\label{subvar}
The level sets of $\kappa$ are smooth subvarieties  of $\RR^{3}$
except for $\kappa^{-1}(0)$ and $\kappa^{-1}(-2)$
each  containing a unique singular point
respectively $(0,0,0)$ and $(2,2,2)$.
\end{lemma}
\proof 
By direct computation one sees that $\kappa$ is a submersion at $(a,b,c)$ unless
 $$ a = b = {c\over c-1} \text{ and }c^{2}-2c = 0.$$
This system has exactly  two solutions,
 namely $(0,0,0) \in \kappa^{-1}(0)$ and  $(2,2,2)\in \kappa^{-1}(-2)$.
The former corresponds to the trivial representation and the latter to  representations such that
$CBA = -I_{2}$.
$\Box$

It follows that the geometric component is 
a smooth subvariety which we shall now describe more fully.
Observe, that  if $ab - a -b \neq 0 $ then 
\begin{eqnarray}\label{c}
c &=& {\kappa - (2 - ab)  \over ab - a -b }
\end{eqnarray}
and we use this formula to determine the topology of the  components
of the relative character variety.

\begin{lemma}[Connected components]
If $t > -2$ then 
$$\kappa^{-1}(t) = \{(a,b,c)\in \RR^3: 2  + abc -ab - bc - ac= t\} $$
consists of three  components each of which is a graph over
one of the components of
$\{ab-a-b \neq 0 \}\subset \RR^2$.
\end{lemma}
\proof
If $ab-a-b = 0$ then  $ab>4$ so $2 - ab - \kappa \neq 0$.
(otherwise there may be a point of indeterminacy in the formula (\ref{c})
for $c$ above and the map $(a,b,c) \mapsto c$ fails to be a diffeomorphism onto
the components of $ab-a-b \neq 0$.) 

The equation $ab - a -b =  0 $ defines a right hyperbola with asymptotes $a=1$ and $b=1$
so that 
$\{ ab - a -b \neq  0 \}$ has three simply connected components, namely:
\begin{enumerate}
\item
 $ab - a -b >0,a,b>1$
 \item
 $ab - a -b <0$
 \item
  $ab - a -b >0,a,b<-1$.
 \end{enumerate}$\Box$

\section{Hyperbolicity and Inequalities}

In this section we establish two inequalities  for functions 
on the geometric component of  $\kappa^{-1}(t >-2)$.
Firstly,  we show that $\tr BA < -2$  so that $BA$ (and by symmetry $CB,AC$) is always  a hyperbolic element.
Secondly, we prove an inequality that shows that a weak  verion of the Collar Lemma holds.

\subsection{Hyperbolicity  of the product $BA$}

We begin by showing that $ab >4$ on the geometric component of the $\kappa>-2$
which implies that $BA$ is always  a hyperbolic element.
Another important  consequence is that no two of $a,b,c$ are simultaneous equal to 2
 on the geometric component,
it follows  from this that the fixed points sets of the restrictions of 
the involutions $I_{a},I_b,I_{c}$ to the geometric component are disjoint.

\begin{lemma}[Product]\label{product}
If  $\kappa>-2$ and $a,b,c>1$ then 
\begin{eqnarray}
\label{fundamental inequality}
ab &>&4,
\end{eqnarray}
so that $$\tr \rho(\alpha\beta) = \tr BA = 2- ab < -2.$$
\end{lemma}
\proof  By the preceding lemma, we  need to  determine a lower bound for
$ab$ over the region $X=\{(a,b)\in \RR^2: a,b>1, ab -a -b >  0\}$,
and this is equivalent to minimizing subject to the  constraints   
$$a>1,\,  b> {a\over a-1}.$$
Thus
$$a b >  \left( a \times {a\over a-1} \right)
=  \left({a^2\over a-1} -4 \right) + 4 
=  {(a-2)^2\over a-1} + 4 \geq 4,$$
since $a-1>0$.
$\Box$

\subsection{Proof of the inequality (\ref{mmc2}) }






From the preceding paragraph the 
elements $\alpha\beta$ and $\beta\gamma$ are hyperbolic 
so that there are positive real numbers 
$\ell_{\alpha\beta}, \ell_{\beta\gamma}$ such that 
$$\tr \rho(\alpha\beta) =   2 - ab  =  -2 \cosh(\ell_{\alpha\beta}/2)),\,
\tr \rho(\beta\gamma) =   2 - ab  =  -2 \cosh(\ell_{\beta\gamma}/2)). $$
It follows  that 
\begin{eqnarray*}
(4 - ab )(4 - bc ) & = & (2  + \tr(\alpha\beta))(2  +  \tr(\beta\gamma))\\
& = & (2 - 2 \cosh(\ell_{\alpha\beta}/2))(2 - 2\cosh(\ell_{\beta\gamma}/2)\\
& = & 16 \sinh^2(\ell_{\alpha\beta}/4) \sinh^2(\ell_{\beta\gamma}/4)
\end{eqnarray*}
We minimize $(ab - 4)(bc - 4)$ over 
the relative character variety
to obtain (\ref{mmcd}).

\begin{lemma}
\label{spCase}
\begin{eqnarray}
(ab - 4)(bc - 4) \geq  4 (\kappa + 2) \label{myCollar}
\end{eqnarray}
\end{lemma}
\proof
Set $ab = h,\, h >4$.
The region $ab - a - b > 0,a>1,b>1$ 
is foliated by arcs  of  the hyperboloids $ab = h$. 
Each leaf meets $ab - a - b = 0$ 
 in exacly 2 points. 
To prove the inequality (\ref{myCollar})it suffices,
to minimize $bc$ along each of the leaves of this foliation.\\
Begin by observing that
$$hc - h - bc - {hc \over b} = \kappa - 2 $$
and consequently,
$$bc = b \left({ \kappa - 2 + h \over h - b - h/b}\right)
= { \kappa - 2 + h \over h/b - 1 - h/b^2}  $$
Minimising $bc$ whilst keeping $h$ constant
 is  equivalent to  maximising the function
$$b \mapsto {h\over b} - 1 - {h\over b^2}
= \left({h \over 4}- 1 \right)  -h\left({1\over b} - {1\over 2}\right)^2 
$$
From the latter expression one sees that 
the function has a unique minimum  when $b = 2$
and the corresponding value of $bc-4$ is
$$4\left( { \kappa - 2 + h \over h - 4} \right)-4
= 4\left( { \kappa +2 \over h - 4}  + 1\right)
= 4\left( { \kappa +2 \over ab - 4}\right)
$$
$\Box$

\section{Dynamics}

We define  involutions of $\fund$ which are the analogues of the ``diagonal exchanges" in the punctured torus
$\alpha,\beta \mapsto \alpha,\beta^{-1}$.
These maps correspond to ``topological reflections" (Section 3 \cite{Nana},
and it is well known that they generate a subgroup of finite index in the mapping class group.
Subsequently, we give explicit formulae for the induced maps
on the character variety and use these to investigate the dynamics.

\subsection{Induced automorphisms of the character variety}
\label{autos}

In this section we show how to calculate the induced map $(a,b,c) \mapsto (a',b',c')$
using cross ratios of  fixed points of parabolics.

Let $f : \pi_1 \ra \pi_1$ be an automorphism. 
The automorphism  $f$ induces a homeomorphisms 
on the Gromov boundary $\partial_\infty \pi_1$ of the fundamental group  which,
for discrete, faithful representations of $\pi_1$, 
induces a $\pi_1$ equivariant homeomorphism $\ff$ on the limit set.
This simply means that if $\gamma\in \pi_1$ acts on $\dhyp$ with a fixed point $\gamma^+$ then
$f(\gamma)$ has a fixed point $f(\gamma)^+$ and
\begin{equation}\label{functorial}
\ff(\gamma^+) = f(\gamma)^+ .
\end{equation}
This relation that allows us to calculate $a',b',c'$ explicitly using the cross ration  on $\RR\cup \{ \infty \}$ as follows.
Recall that $\alpha, \beta,\gamma$ 
fix $0,1,\infty$ respectively,
so we set
$$\alpha^+ = 0,\, \beta^+ = 1,\, \gamma^+ = \infty.$$
We shall write each of  $a,b,c$ as (a function of) a cross ratio of these three  points plus another point $x\in \RR$.
For $x\neq 0,1,\infty$, one has the following identity 
\begin{equation}\label{identity}
x = [x, 1, 0, \infty] = { x - 0  \over x - \infty }.{ 1- \infty  \over 1 - 0 } =  
{  x - \alpha^+  \over x  - \gamma^+ }.{ \beta^+ - \gamma^+ \  \over \beta^+ - \alpha^+ }.
\end{equation}
Now observe that
\begin{eqnarray*}
a   =  {1 \over  A(\infty)}  &=& {1 \over \alpha(\gamma^+) } \\
c  =   C(0) &=&  -\gamma(\alpha^+)\\
b  = {1\over B(\infty) - 1 }  &=& {1 \over \beta(\gamma^+) - 1}.
\end{eqnarray*}

The following lemma is an immediate consequence of these three equations and (\ref{functorial}),(\ref{identity}).

\begin{lemma}
Let $f : \pi_1 \ra \pi_1$ be an automorphism.
Then $f$ induces a map
$$
\begin{pmatrix}
a\\
b\\
c
\end{pmatrix} \mapsto 
\begin{pmatrix}
a'\\
b'\\
c'
\end{pmatrix}
$$
on the set of $(a,b,c)$ such that the corresponding representation is discrete faithful.
The values of $a',b,,c'$ are given by:
\begin{eqnarray*}
a'& = &{ 1 \over F(f(\alpha)(f(\gamma)^+))} \\
b'& = &{1 \over F(f(\beta)(f(\gamma)^+))) -1) }\\
c'& = &- F(f(\gamma)(f(\alpha)^+)))
\end{eqnarray*}
where the  function $F : \RR \setminus \{0,1 \}$ is defined by
$$ F(x) = 
{ x  - f(\alpha)^+  \over x  - f(\gamma)^+ }
.{ f(\beta)^+ - f(\gamma)^+   \over f(\beta)^+ - f(\alpha)^+ }.$$
\end{lemma}
We apply the lemma to determine the  map induced by the involution $\inv_b$.

\begin{coro}
The  action of the map $\inv_{b}$
induced by the automorphism 
$\phi_\beta$ is
$$
\inv_{b} : \begin{pmatrix}
a\\
b\\
c
\end{pmatrix}
\mapsto 
\begin{pmatrix}
a(b-1)\\
b/(b-1)\\
 c(b-1).
\end{pmatrix}$$
\end{coro}

\proof
The morphism  $\phi_\beta$ is defined by
\begin{eqnarray*}
\alpha &\mapsto& \alpha^{-1}\\
\beta &\mapsto& \alpha^{-1}\beta^{-1}\alpha\\
\gamma &\mapsto& (\beta\alpha)^{-1}\gamma^{-1}(\beta\alpha).
\end{eqnarray*}

\begin{center}
\includegraphics[scale=.4]{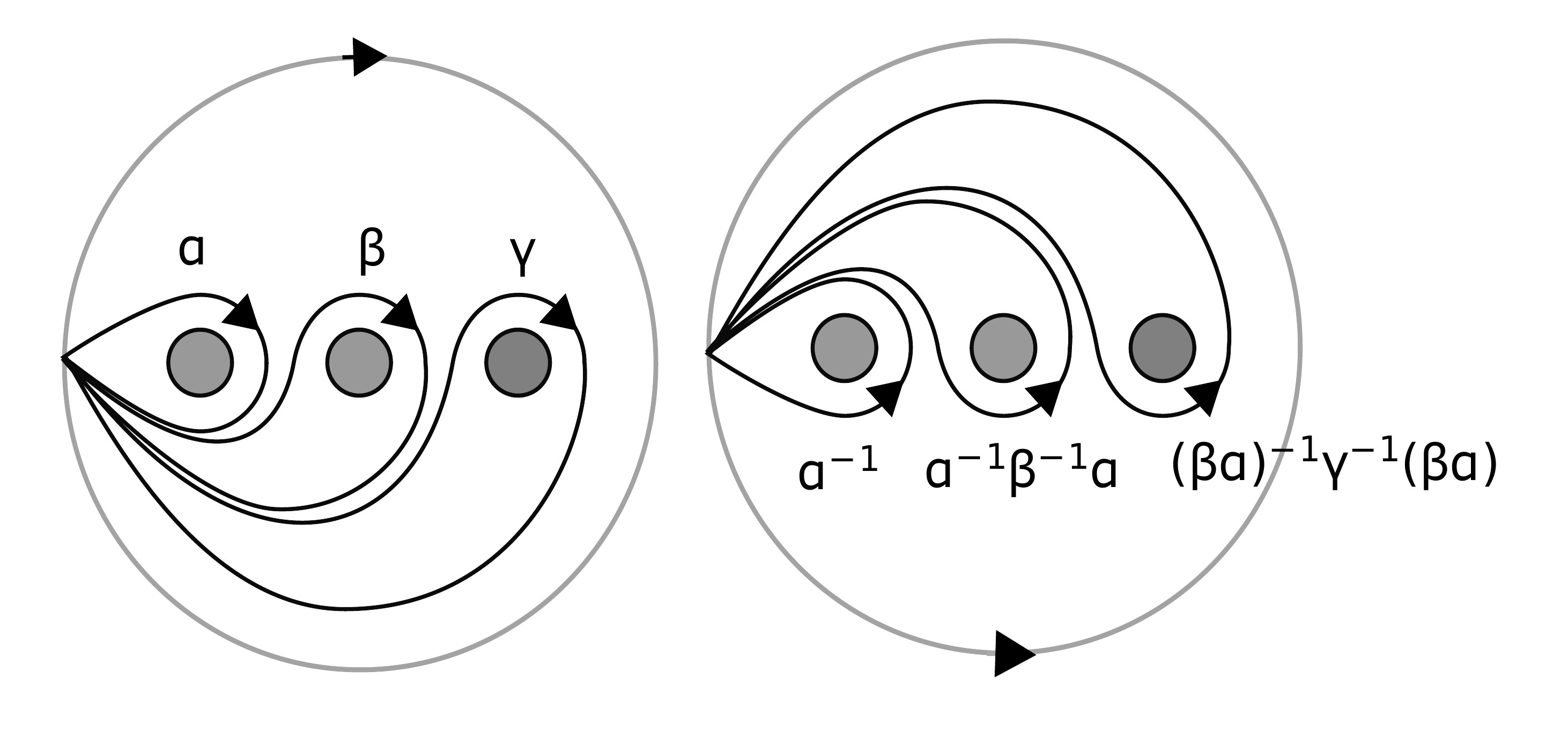} 
\end{center}

We use these equations to determine fixed points as follows
\begin{eqnarray*}
 f(\alpha)^+ &=& (\alpha^{-1})^{+} = \alpha^+ = 0 \\
 f(\beta)^+ &=& \alpha^{-1}(\beta^+) =   \alpha^{-1}(1) = {1 \over 1-a } \\
 f(\gamma)^+ &=& (\beta\alpha)^{-1}(\gamma^+) = (\beta\alpha)^{-1}(\infty) =  {1-b \over  ab - a - b  }\\
\end{eqnarray*}
One also has, by direct computation,
$$f(\gamma)(f(\alpha)^+) = f(\gamma)(0) = \frac{-b^{2} c+2 b c-c}{a b^{2} c-2 a b c+a c-b^{2} c+b c-1}.$$
It is then straightforward to check,
 either by hand or using a computer algebra package,
 that the morphisms induced is $I_b$. $\Box$

\subsection{Fixed point sets of induced automorphisms}



Each of the three involutions $\inv_a,\inv_b,\inv_c$, has a non empty fixed point set.
The dynamics of the group generated by these  three involutions is 
determined by the geometric configuration of these fixed point sets

\begin{lemma}.
\begin{enumerate}
\item Viewed as a self map  of $\RR^3 \setminus \{ (b -1)\}$, the  fixed point set of the map $\inv_{b}$ 
consists of the plane $\Pi_{b} := \{b=2\}$ and $(0,0,0)$.
\item 
The plane $\Pi_{a}$ (resp. $\Pi_{b},\Pi_{c}$)
meets the level sets $\kappa^{-1}(t>-2),a,b,c>1$ in a hyperbola.
The three hyperbolae obtained in this way are disjoint.
\item
The three planes  meet the level set $\kappa^{-1}(t=-2),\,a,b,c>1$  in the lines:
$$\{a=2\}, \{b=2\}, \{c=2\}.$$
These lines are concurrent at the singular point of $\kappa$,  $(2,2,2)$.

\end{enumerate}
\end{lemma}

\proof 
For the first point, observe that if $(a,b,c)$ is a point of the fixed point set of $I_b$ then one has
$$ b = {b \over b - 1},\,  a = a(b-1),\, c = c(b-1).$$
Clearly, if $b \neq 0 $ then $b -1 = 1$ and if $b=0$ then $a = -a, c = -c$ so that $a = b = c = 0$.

For the second point, it is convenient given  to consider $I_c$.
If $c = 2$ then $\kappa = ab - 2a - 2b + 2 = (a-2)(b-2) - 2  $ and this is 
the equation of a right hyperbola in the plane $c  = 2$ provided $\kappa + 2 > 0$.

Finally, if $\kappa + 2 > 0$ the hyperbola $c = 2$ is asymptotic to the lines $ a = 2, b = 2$
and so the hyperbolae are disjoint as required. 
If $\kappa  =  - 2$ the intersection of $c = 2$ with the level set is the pair of lines $ a = 2, b = 2$
and the third point of the lemma follows immediately.
$\Box$

\subsubsection{The action is proper $\kappa>-2$}

We apply the preceding lemma to show that the 
group of automorphisms acts properly discontinuously  
and to  determine  a fundamental domain for this action.

\begin{thm}[Proper]
The three involutions generate a group $\Gamma$
isomorphic to $\ZZ/2\ZZ * \ZZ/2\ZZ * \ZZ/2\ZZ$
which acts properly on the geometric component.
Furthermore a $\Gamma$-fundamental domain is
$$\Delta = \{(a,b,c): \min(a,b,c)>2\}.$$
\end{thm}

\begin{center}
\includegraphics[scale=.6]{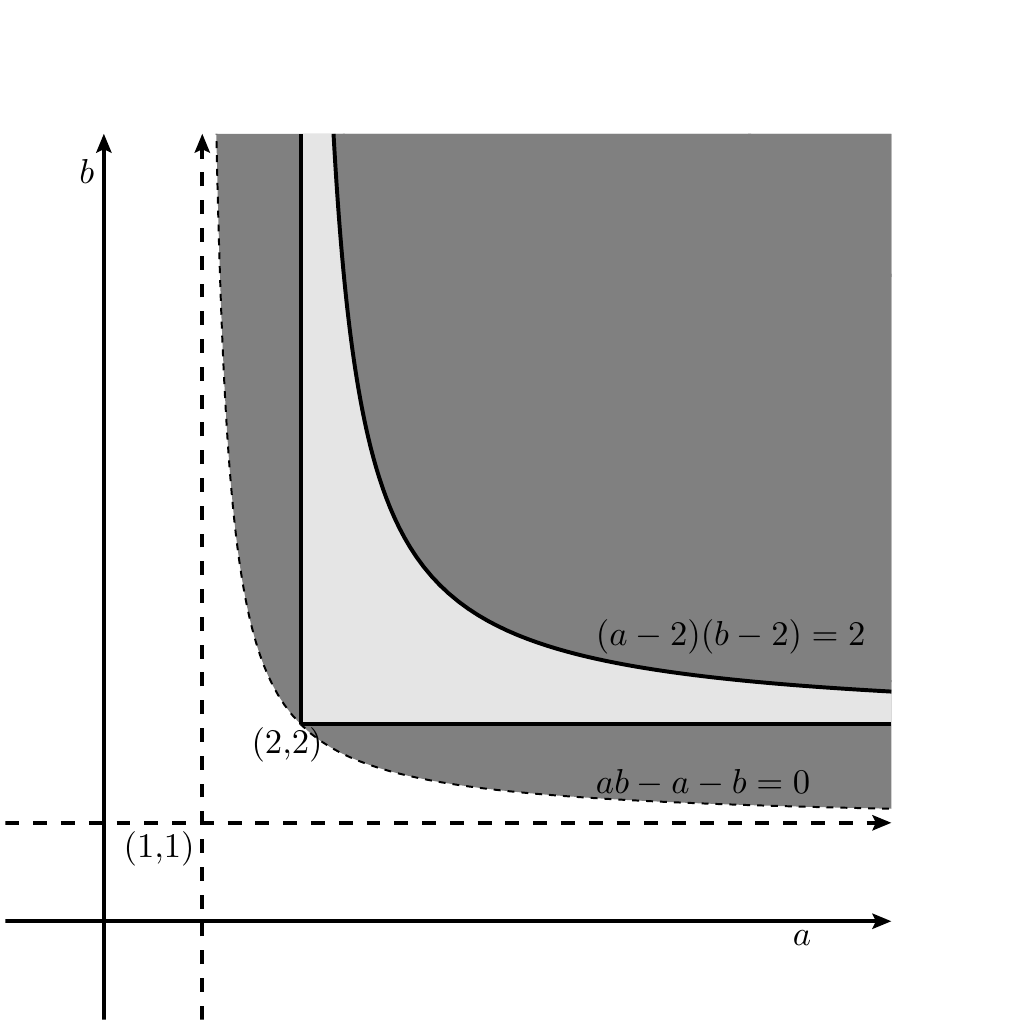} 
\end{center}

\proof To prove that $\Delta$ is a fundamental domain we have to show 
the following:
\begin{itemize}
\item $\Delta,\inv_{a}(\Delta),\inv_{b}(\Delta),\inv_{c}(\Delta)$ are disjoint.
\item if  $x\in\kappa^{-1}(t>-2)$ then there is a word in $w\in\Gamma$ such that 
$w(x)\in \overline{\Delta}$
\end{itemize}
For the first of these points, we note that 
$$a>2 \Rightarrow {a\over a-1}<2$$
 so $\Delta \cap \inv_{a}(\Delta) = \emptyset$
 (and by symmetry  $\Delta \cap \inv_{b}(\Delta), \Delta \cap \inv_{c}(\Delta)  = \emptyset$.)
 
For the second point suppose $(a,b,c) \notin \Delta$
and, imitating Wolpert \cite{??} for the punctured torus,
we choose  an ``energy" function $E$ on the level set.
An \textit{energy function} is a proper function   
such that  for any   $(x,y,z)$ not in the closure of the  fundamental domain $\overline{\Delta}$,
there one of the involution $\inv_{a},\inv_{b},\inv_{c}$
which   decreases its value, that is
$$ \min (
E\circ\inv_{a}(x,y,z),E\circ\inv_{b}(x,y,z),E\circ\inv_{c}(x,y,z)
) <E(x,y,z)$$
Note that 
$$(\overline{\Delta})^c = \{ (a,b,c): \min(a,b,c) <2 \}.$$ 
The function
$$E:(a,b,c) \mapsto abc,$$
satisfies this second condition since
\begin{eqnarray}
 E\circ \inv_{a}(a,b,c) = abc(a-1) < E(a,b,c)
\end{eqnarray}
if $1<a<2$.
The function $E$ is clearly continuous so 
 it is proper if the preimage of every  bounded set is bounded.
By Lemma \ref{product}, if $\kappa>-2$ then,
  $ab,bc,ac>4$
it follows immediately that, for any $K>0$,
$$E(a,b,c)\leq K \Rightarrow  \max\{a,b,c\} \leq K/4.$$
$\Box$

\section{Symplectic volume}

Wolpert has shown that the symplectic form takes the form
$d \ell_\alpha \wedge d \tau_\alpha$  where $\ell_\alpha$ is the length as before
and $\tau_\alpha$ is a Fenchel-Nielsen twist parameter.
We derive an expression  for this in terms of $a,b,c$ and use it to calculate the volume.

\subsubsection{Fenchel Nielsen coordinates}

From the preceding discussion we have
$$ 2 \cosh(\frac12 \ell_\alpha)  = ab - 2 = - \tr \alpha$$
Then

\begin{eqnarray*}
d  (2 \cosh(\frac12 \ell_\alpha)) &=& \sinh(\frac12 \ell_\alpha)  d \ell_\alpha \\
&=&\frac12 ((\tr \alpha)^2 - 4)^\frac12   d \ell_\alpha) \\
 &=&\frac12 \Delta(ab) d \ell_\alpha
\end{eqnarray*}
 where  $\Delta(ab)  := \sqrt{(ab - 2)^2 - 4} =  2\sinh(\frac12 \ell_\alpha) $.
 
 So  $d \ell_\alpha = \frac{2}{\Delta} d(ab)$ depends only on $ab$
 and we will exploit this by writing $d \tau_\alpha$ as a sum
 where one term is $F(ab)\mathrm{d} (ab)$ for some function $F:\RR \rightarrow \RR$.

The twist parameter is given by the signed distance between 
the image of a pair of reference points under the projection 
to the axis of $\alpha$. 
We choose the  fixed points of $A,C$
(respectively $0$ and $\infty$ )
 as reference points  and then
$$\tau_\alpha = 
 \log \left(\frac{\alpha^+ -0 }{\alpha^- - 0 }.\frac{\alpha^+ - \infty }{\alpha^- -\infty } \right )
 = \log(\frac{\alpha^+}{\alpha^-} )$$
Now, by direct calculation,  we express this as a function of $a,b$
$$\alpha^\pm = \frac{1+ b  - ab - (1-b) \pm \Delta}{b+a - ab}$$
so that 
$$d \tau_\alpha = d \log ( 2b - ab +  \Delta) - d \log (2b - ab - \Delta)  = 
\frac{2\Delta  (2\mathrm{d}\,b + F(ab) \mathrm{d}(ab)) }{(2b - ab)^2  -   \Delta^2 }  , $$
for some function $F$. Thus
$$d\ell_\alpha \wedge d \tau_\alpha =
 \frac{\mathrm{d}\,b \wedge \mathrm{d}(ab) }{b(ab - a - b) } = {db \wedge da \over ab - b - a}.$$

\subsubsection{Volume}
The expression for the volume form in these coordinate is just:
$$\omega_{WP} = {da \wedge db \over ab - b - a}= {db \wedge dc \over bc - b - c} = {dc \wedge da \over ac - a - c}.$$

$$\int_2^\infty \int_2^{\frac{2-2 b-k}{2-b} }
 {da \over ab - a - b } {db}.$$
 For the particular case $k =2$ this becomes
 the Fundamental domain is a region between two hyperboloids
$(a-2)(b-2)  = 0$ and $(a-2)(b-2)  = 4$
so the integral becomes
 $$\int_2^\infty \int_2^{a = f(b) }
 {da \over (a-1)(b-1) -1  } {db},$$
 where $(f(b)-2)(b-2) = 4$.
 Using the change of variable $u = a-2,\,\, b = v-2$ 
  \begin{eqnarray*}
  \int_0^\infty \int_0^{4/v}
 {du \over (u+1)(v+1) -1  } {dv} 
 &=& \int_0^\infty  {2\over v + 1} \log( {v+2 \over v}) {dv} \\
 &=&2 \int_0^\infty  {\log( v+2)  \over v + 1} {dv}  - 2\int_0^\infty  {\log (v) \over v + 1} {dv} \\
   &=& 2[-Li_2(-t-1)-Li_2(-t)-\log(t) \log(t+1)]_0^\infty   \\
  &=& {2\pi^2 \over 3} -  {\pi^2 \over 6}   =  {\pi^2 \over 2}
 \end{eqnarray*}
  So the volume of the quotient of the 
 component of the character variety of the 4 
 by the group generated by the involutions is $\pi^2/2$.
 By considering the action on the abelianisation of the fundamental group
 one sees that the mapping class group is a finite subgroup of order $4$
 in this latter group so that volume of the moduli space is $2\pi^2$
 which agrees with the value of $V_0(0,0,0, i\theta ) = {1 \over 2}(4 \pi^2  - \theta^2)$
 when $\theta = 0$.

\section{Hyperbolicity implies Fibonacci Growth}

The free product $(\ZZ/2\ZZ) * (\ZZ/3\ZZ)$ 
has a Bass-Serre tree for which the 
edge stabiliser  are isomorphic to $\ZZ/2\ZZ$
and the vertex stabiliser  are isomorphic to $\ZZ/3\ZZ$.
One identifies the Markoff tree with the Bass-Serre tree
by identifying $PSL(2,\ZZ)$ with the (smallest) group of automorphisms 
of the Markoff cubic that is transitive on ordered triples of integers $(x,y,z)$ 
that are  solutions. 
The vertices of the Markoff tree are the unordered triples  $\{ x,y,z \}$ 
and the edges unordered pairs $\{x,y\}$.
One can embed the Markoff tree,
viewed as the Basse-Serre tree of $PSL(2,\ZZ)$,
in the Poincar\'e disc. 
The complement consists of countably many \textit{complementary regions}
three of which meet at every vertex $\{ x,y,z \}$ 
and  one of these numbers can be associated  to each region in a consistent manner 
(see Bowditch \cite{bow} or  \cite{tanAdv} for details).

%

Using this construction, Bowditch defined   \textit{Fibonacci growth} 
for the $PSL(2,\ZZ)$ orbit of solutions.
Let $z_0 \geq y_0 \geq x_0$ be a vertex and
consider the  binary subtree 
which is union of the edge $e = \{x_0,y_0 \}$ and  the component of $ T \setminus  \text{int} (e)$  containing the vertex $ \{x_0,y_0, z_0\}$.
We think of this subtree as starting at the edge  $e = \{x_0,y_0 \}$.
One introduces a comparison function defined recursively 
on (complementary regions of an embedding in the disc of) 
such a binary subtree.
If  $\{X,Y\}$ is the first edge one sets 
$$F_e(X) = \log(x_0),  F_e(Y) = \log(y_0)$$
If the function is defined at the edge 
$\{X,Y\}$ and $\{X,Y,Z\}$ is a vertex 
such that $F_e(Z)$ is as yet undefined  then 
$$ F_e(Z) = F_e(X) + F_e(Y).$$
A function $f$ on the vertices admits  \textit{upper and lower Fibonacci bounds}
if there exists $K^-,K^+ >0>0$ such that for every vertex $v$ one has
\begin{equation}\label{fib}
	K^- F_e(v) \leq f(v) \leq K^+ F_e(v).
\end{equation}
An important observation of Bowditch is that 
if there exists a constant $\delta>0$ 
such that a real valued  function $f$ defined on the vertices of the tree satisfies 
\begin{equation}\label{bowditch}
	f(z) \geq f(x)  + f(y) - \delta 
\end{equation}
at each vertex  $z \geq y \geq x$ in a binary sub tree starting at $z_0 \geq y_0 \geq x_0$
then 
\begin{equation}
	f(z) \geq  (m - \delta) F_e(z)  + \delta,
\end{equation}
where $m = \min(x_0,y_0,z_0)$. 
In particular it admits a lower Fibonacci bound.

\begin{lemma}
The function $f :  (a,b,c) \mapsto \log(ab)$
satisfies the inequality (\ref{bowditch})  above with $\delta \leq \log(4)$.
\end{lemma}

\proof 
It is convenient to adopt the following notation:
for a function   $(a,b,c) \mapsto f(a,b,c),\, 
I_b^*(f(a,b,c)) := f\circ I_b(a,b,c)$.
Under this operation the functions 
\begin{eqnarray*}
(a,b,c) &\mapsto& ab,\\
(a,b,c) &\mapsto& bc,\\
 (a,b,c) &\mapsto& ac
\end{eqnarray*}
transform as follows:
\begin{itemize}
	\item $I_b^*(ab) = ab , I_b^*(bc) = bc$ that is these are invariant functions 
	\item  $I_b^*(ac) =  ac(b-1)^2$.
\end{itemize}
Consider the ratio
$$ {I_b^*(ac) \over I_b^*(ab) I_b^*(bc)} 
= { ac(b-1)^2 \over ab^2c}  
= { (b-1)^2 \over b^2} = \left(1- {1\over b} \right)^2.$$
Taking logs:
\begin{eqnarray*}
 \log(I_b^*(ac))  
&=& \log( I_b^*(ab))  + \log( I_b^*(bc)) + 2\log (1- {1/b})\\
&=& \log(ab)  + \log(bc) - 2\log \left({b \over b - 1} \right) 
\end{eqnarray*}
By hypothesis  $b  > 2$ so that $2\log (b/  (b - 1) ) < 2 \log(2) = \log(4)$ 
and it follows that this system satisfies Bowditch's condition (\ref{bowditch}) if
$$\delta = \log(4) \leq  \min(\log(ab),\log(bc)).$$
This is equivalent to 
$$4 \leq  \min(ab,bc),$$
which is immediate from Lemma \ref{product}. $\Box$

\section{A hyperbolic metric with holonomy $\rho$}

Finally we show that each representations 
is the holonomy of a singular hyperbolic metric on the 
3 punctured sphere. 
This is done by constructing a convex 
polygon $P_\theta$  in $\HH$
such that  
\begin{enumerate}
\item each of $A,B,C$ act as side pairing transformations
\item the quotient space obtained by identifying sides 
is a 3 punctured sphere.
\end{enumerate}
The $P_\theta$  has 3 ideal  vertices,
namely the fixed points of $A,B,C$,
and 3 finite  vertices
which are the fixed points of the elliptic elements
$CBA$, $BAC$ and $ACB$ respectively.

\subsection{The fixed point of $CBA$. }


\begin{lemma}\label{sign}
If $a + b -ab <0,\, b>2, 2 >\kappa >-2$ then 
$CBA$ is elliptic and its 
 unique fixed point $z\in \HH$ 
satisfies 
$$ \text{Re} (z_\pm)  
=  { \kappa  - 2 + 2b  \over a + b -ab  } < 0 $$
\end{lemma}
\proof
One computes  $CBA$
$$CBA = \begin{pmatrix} 
1 + b - ab  -ca -cb  + abc &   -b -c + bc\\
a + b - ab   &  1-b
\end{pmatrix}$$
and applies the formula (\ref{fpsign})$\Box$

\subsection{Proof of existence.}

\begin{thm}\label{metric}
Provided $\kappa(a,b,c) > -2$ the representation 
$\rho$ associated to the triple $(a,b,c)$ is the holonomy of a hyperbolic metric   with  cone angle $\theta$ such that 
$$2\cos \theta/2 = \kappa.$$
\end{thm}
\proof
Let $z\in \HH$ be the fixed point of $CBA$. 
Begin by observing that $C^{-1}$ conjugates $CBA$ to $C^{-1}(CBA)C = BAC$
and
$BA$ conjugates $CBA$ to $BA(CBA)(BA)^{-1} = BAC$.
Thus we need to check that the ``fundamental polygon"
with the following 6 vertices 
$z, 0, C^{-1}(z), (BA)^{-1}(0),1,0$ is convex.

This follows from Lemma \ref{sign}.
$\Box$

%


\begin{thebibliography}{99} 

     \bibitem{bow} 
   B.H. ~Bowditch
   Markoff triples and quasifuchsian groups
   Proc. London Math. Soc. (1998) 77 (3): 697-736.
   

  \bibitem{DoNo} 
Norman ~Do, Paul  ~Norbury
Weil-Petersson volumes and cone surfaces
Geometriae Dedicata 141 (2009), 93-107.



  \bibitem{Tan} 
I. Kim, J. Kim, S.P. Tan
McShane's Identity in Rank One Symmetric Spaces, Math. Proceedings of the Cambridge Philosophical Society, 157 (2014), 113--137.


    \bibitem{MaPaTa} 
  S. ~Maloni, F. ~Palesi,  S.P. ~Tan, 
  On the character variety of the four-holed sphere,
   eprint arXiv:1304.5770, 

  \bibitem{maskit} 
B. ~Maskit
 "Matrices for Fenchel-Nielsen coordinates.."
  Annales Academiae Scientiarum Fennicae. Mathematica 26.2 (2001): 267-304
  
    \bibitem{NaNa} 
T. ~Nakanishi  M. ~Naatanen
Areas of two-dimensional moduli spaces
Proc. Amer. Math. Soc. 129 (2001), 3241-3252
  

  \bibitem{penner} 
R.C. ~Penner,  
The decorated Teichmueller space of punctured surfaces. 
Comm. Math. Phys. 113 (1987), no. 2
   
     \bibitem{SchTra} 
  
Georg ~Schumacher,  Stefano ~Trapani. 
Weil-Petersson geometry for families of hyperbolic conical Riemann surfaces. Michigan Math. J. 60 (2011), no. 1, 3--33. 


  \bibitem{tanzhang} 
S.P. ~Tan, Y.  ~Wong, and Y. ~Zhang
Generalizations of McShane's identity to hyperbolic cone-surfaces
J. Differential Geom.
Volume 72, Number 1 (2006), 73-112.

  \bibitem{tanAdv} 
S.P. ~Tan, Y.  ~Wong, and Y. ~Zhang
Generalized Markoff maps and McShane's identity. 
Advances in Mathematics 217:2 (2008)



  \bibitem{McLab} 
   F. ~Labourie,G.  ~McShane
Cross ratios and identities for higher Teichmüller-Thurston theory
   Duke Math. J.
Volume 149, Number 2 (2009), 279-345.

  \bibitem{wolp} 
S. ~Wolpert
On the Kaehler form of the moduli space of once punctured tori.
Commentarii mathematici Helvetici (1983)
Volume: 58, page 246-256

\end{thebibliography}
\end{document}